\documentclass[12pt,draft]{amsart}

\newtheorem{theorem}{Theorem}[section] 

\newtheorem{thm}[theorem]{Theorem}
 
\newtheorem{proposition}[theorem]{Proposition} 
 
\newtheorem{corollary}[theorem]{Corollary}

\theoremstyle{definition}

\newtheorem{conj}[theorem]{Conjecture} 
\newtheorem{prob}[theorem]{Problem} 
\newtheorem{problem}[theorem]{Problem}

\theoremstyle{remark}

\DeclareMathOperator{\cf}{cf}
\DeclareMathOperator{\pcf}{pcf}
\DeclareMathOperator{\pp}{pp}
\DeclareMathOperator{\bd}{bd}

\newcommand{\aaa}{{\mathfrak a}} %%% a di pcf a 
\newcommand{\pcfa}{\pcf \aaa}

\title[decomposable ultrafilters and possible cofinalities]
{A connection between decomposability of ultrafilters and possible cofinalities}

\author[]{Paolo Lipparini} 
\address{Dipartimento di Matematica\\
II Universit\`a di Roma (Tor Vergata)\\
Valle di Lacrime Scientifiche\\
I-00133 ROME 
ITALY
}
\email{lipparin@axp.mat.uniroma2.it}
\urladdr{http://www.mat.uniroma2.it/\textasciitilde lipparin}
\thanks{The author has received support from MPI and GNSAGA
} 


\begin{document} 

\begin{abstract} 
We introduce the decomposability spectrum 
$K_D=\{\lambda \geq \omega| D \text{ is }
\lambda\text{-decomposable}\}$ of an ultrafilter $D$, and show that Shelah's 
$\pcf$ theory influences the possible values $K_D$ can take.

For example, we show that  if 
$\aaa$ is a set of regular cardinals, 
$\mu \in \pcfa$, 
the ultrafilter $D$ is $|\aaa |^+$-complete and $K_D \subseteq  \aaa$, then 
$\mu \in K_D$.

As a consequence, we show that if $ \lambda $ is singular and
for some $ \lambda' < \lambda $ $K_D$ contains all regular 
cardinals in $ [\lambda' , \lambda)$
then: 

(a) if $\cf \lambda = \omega $ then
 either $ \lambda \in K_D$, or $ \lambda ^+ \in K_D$; and

(b) if $D$ is  $(\cf \lambda)^+$-complete then $ \lambda ^+ \in K_D$,
and $\pp ( \lambda )= \lambda ^+$.
\end{abstract}

\maketitle 

\section{Introduction} \label{intro}

An ultrafilter $D$ over $I$ is {\it $\lambda$-decomposable} 
if and only if there is a partition of $I$ into $\lambda $
sets such that the union of any $<\lambda $ sets of the partition
never belongs to $D$. 
%togliere in quello per JSL??? Citare applicazioni?
In other words, $D$ is $\lambda $-decomposable 
if and only if some quotient of $D$ is uniform over $\lambda $;
to be more precise, $D$ is $\lambda $-decomposable 
if and only if there exists some ultrafilter  $D'$
uniform over $ \lambda $, and  $D'\leq D$ in the 
{\em Rudin-Keisler order}.
If $D$ is an ultrafilter, define the {\it decomposability spectrum} 
$K_D$  of $D$ by $K_D=\{\lambda \geq \omega| D \text{ is }
\lambda\text{-decomposable}\}$. 
In this note we address the following question: which are the possible
values $K_D$ can take?

It is well known that $K_D$ is closed under taking cofinalities and
regular predecessors; more explicitly, if $ \kappa $ is regular and $ \kappa^+ \in K_D$ 
then $ \kappa \in K_D$; and, if $ \kappa \in K_D$ is singular, then 
$\cf \kappa \in K_D$. Further constraints on $K_D$ are given in \cite{L1}. In this
 note we show that, under a completeness assumption on $D$,
$K_D$ is closed under Shelah's $\pcf$ operation of taking possible
cofinalities of reduced products.

We now give a few examples of possible values for $K_D$.
The possibility that $K_D$ is an interval can always occur:
if $D$ is uniform over $ \lambda $ and $( \omega , \lambda )$-regular  then
$K_D= [\omega,\lambda]$. By results from \cite{D} (see also \cite{L1}),
if there is no inner model with a measurable cardinal then
 $K_D$  is always an interval with $ \omega  $ as the 
inferior extreme.

If there is a measurable cardinal $\mu$, then there is
a $\mu$-complete ultrafilter $D$ over $\mu$, and this implies 
$K_D=\{\mu\} $. 
Conversely, if $|K_D|=1$, say $K_D=\{\mu\} $,
then $\mu $ is either $\omega $ or a measurable cardinal.
Starting from a measurable cardinal, we can create 
ultrafilters with gaps in their
decomposability spectra: as a very elementary example, if 
$D'$ is non principal over $ \omega  $, and $D$ is as above, then
 $K_{D\times D'}=\{\omega, \mu\} $.  
Of course, we get more interesting
examples by destroying the measurability of
$\mu$; for example, by Prikry forcing \cite{P}, one can turn the cofinality of $\mu$ to 
$ \omega  $, thus having an ultrafilter for which
$K_D=\{\omega ,\mu\} $ and $ {\rm cf}\mu=\omega $.
%%% see???
By a more elaborate forcing, one can obtain
$K_D=\{\omega ,\mu\} $, for some  strongly
inaccessible and not weakly compact $\mu$ \cite{Shr}.

%mettere? pi- enfasi sul fatto che il mio problema generalizza quello di Si!
The question of which values $K_D$ may
take is intriguing
 even in case $|K_D|=2$. This particular question originated from
\cite{Si}, and
not everything is known yet about those $\mu$ for which we can have 
$K_D=\{\omega ,\mu\} $. Some restrictions on $\mu$ are
listed in \cite
%Section ???
{L1}.
 Seemingly, the possibility $K_D=\{\lambda  ,\mu\} $,
where $ \lambda,\mu$ are both $>\omega $ has never been investigated,
apart from trivial cases \cite
%Section ???
{L1}.

However, in the present note we deal with the case when $K_D$ is infinite.
In this case things are even more involved,
but for a simple reason:
suppose that 
$ \lambda $ is a limit cardinal, 
$ (\lambda_\alpha)_{\alpha\in\cf \lambda}$
 is an ascending
sequence of cardinals unbounded in $\lambda $,
and $ \lambda_ \alpha \in K_D$, for $\alpha\in\cf \lambda$. 
Without loss of generality, $D$ is uniform, say over 
$ \mu  $; this implies that 
$ \mu \in K_D$, and $ \mu  $ is necessarily larger than all the 
$\lambda_\alpha$'s.
The main complication arises from the fact
the sequence 
$(\lambda_\alpha)_{\alpha \in \cf \lambda  }$
conveys unclear and confused information  about
the possible values  $ \mu $ 
can take. 

For example, if we are in the above situation, and
$D$ is $ ( \omega ,\cf \lambda )  $-regular, then $D$ 
is necessarily $ \lambda $-decomposable \cite
%Section ???
{L1}; on the other hand, 
%esempio migliore? there are many cases 
%in which we reach the conclusion that??
if $ \lambda $ is limit,  
 $D$ is $ (\cf \lambda )^+$-complete, and  
% $ \lambda ^{ \cf \lambda } = \lambda ^+$
$\pp( \lambda )= \lambda ^+$
then $D$ is $ \lambda ^+$-decomposable 
(Proposition \ref{solovaygenpp} (d) $ \Rightarrow $ (b)). 
In \cite{L1} we asked whether the above possibilities 
are the only ones which can occur,
namely, we asked the following problem.

\begin{problem}\label{main}
Suppose that $ \lambda $ is a limit cardinal
and that there are arbitrarily large cardinals  $ \lambda ' < \lambda $ 
such that the ultrafilter $D$ is $ \lambda  $-decomposable. 

Is it true that $D$ is either $ \lambda $-decomposable  or
$ \lambda ^+$-decomposable? 
\end{problem}

In other words, Problem \ref{main}
asks whether $K_D$ satisfies the following 
closure property: for every  $X \subseteq K_D$
either $\sup X \in K_D$ or 
$(\sup X)^+ \in K_D$.

In this paper we provide new evidence that
Problem \ref{main} has an affirmative answer
in the great majority of cases when $ \lambda $ is singular. 
Our results, together with $\pcf$ theory,
suggest that a negative solution of Problem 
\ref{main}  for some singular $ \lambda $
could occur only in very special situations, and would probably exhibit a very complicated structure for $K_D$.
%See? O: we shall not give the very cumbersome details for this???
In turn, the existence of an ultrafilter with certain specified values of $K_D$
puts constraints on $ \pcf$ theory and, in some cases, implies that
$\pcf$ theory has the simplest possible structure (Theorem \ref{limit} and Corollary \ref{solovaygenpp}). 

The starting result of the present paper is Theorem \ref{ldec} in the next section.
It shows that $K_D$ satisfies the following closure property:  if
$(\lambda _j)_{j\in J}$ are regular cardinals, $D$
is a $|J|^+$-complete ultrafilter, and $ \lambda_j \in K_D$ for every 
$j \in J$, then $ \mu \in K_D$, for every 
$ \mu$ which can be obtained as the  cofinality of 
$\prod_E \lambda _j$ for a suitable ultrafilter $E$
over $J$.   
A particular case of Theorem \ref{ldec} is rephrased using the
terminology of $\pcf$ theory 
in Corollary \ref{pcfl}. Then we refine the
above mentioned closure property of $K_D$ 
and show that, in many cases, if 
$ \lambda_j \in K_D$ for every 
$j \in J$, then $ \lambda ^+ \in K_D$,
where $ \lambda = \sup _{j \in J} \lambda _j $.
In particular, in Corollary \ref{cfomega}, we show that Problem 
\ref{main} has an affirmative answer in the particular case
in which $X$ is an interval of regular cardinals whose supremum is
singular of cofinality $ \omega $.   
In Section \ref{Sol} we give an alternative proof
for some steps in R. Solovay's result that GCH (the Generalized Continuum
Hypothesis) holds at strong limit singular cardinals
above a strongly compact cardinal. The results in
Section \ref{Sol} do not depend on Section \ref{pcfth}, and their proofs
do not use $\pcf$ theory.
At the end, we state some further problems in Section \ref{prob}. 
See \cite{L1, She} for unexplained notions.  
%definition di $\kappa $-complete. (si deduce da $K_D$!)
%$f$ $ \lambda $-decomposition
%definition di $(\lambda,\lambda )$-regular coi sets [basta $\lambda $-d.i. ?]
%equivalent to $\lambda $-dec se $\lambda $ regular

If $E$ is an ultrafilter over $J$, and $ (\lambda _j)_{j \in J}$
are cardinals, we shall write $ \cf \prod_E \lambda _i$
to denote the cofinality of the linear order
$\prod_E \langle  \lambda_i, \leq\rangle$. 

If $ \lambda $ is a limit cardinal, the locution ``the ultrafilter $D$ is
$ \kappa $-\hspace{0 pt}decomposable for all sufficiently large regular $ \kappa < \lambda $''
means that there exists 
$ \lambda' < \lambda $ such that 
 $D$ 
is $ \kappa $-decomposable
for every regular $ \kappa  $ with $ \lambda  ' \leq \kappa < \lambda $.
%It is a major problem whether in the results of the present paper
%the above assumption can be weakened to %%poi, se c'e' problema di spazio,
%toglierlo dall'ultima sezione!
We shall sometimes consider also 
the weaker condition
``there are arbitrarily large (regular) cardinals  $ \kappa  < \lambda $ 
such that the ultrafilter $D$ is $ \kappa  $-decomposable'',
which means that  for every 
$ \lambda ' < \lambda $ there is some (regular) $ \kappa $
such that $ \lambda ' \leq \kappa < \lambda $, and $D$ is
$ \kappa $-decomposable.
% la prima condizione sembra piu' forte, 
%almeno nella facilita di ottenere dimostrazioni!

We shall make use of the following classical result.

\begin{thm} \label{cckp}
\cite[Theorem 1]{CC} \cite[Theorem 2.1]{KP}
If the ultrafilter $D$ is uniform over $\lambda^+$ then $D$ is either
$\cf \lambda $-decomposable,  or
$(\lambda',\lambda^+ )$-regular for some regular $\lambda'\leq\lambda $.
\end{thm}

%serve? E' equivalente a come  lo mettono loro per cose da dire sopra
See \cite{L1} for further remarks about Theorem \ref{cckp},
as well as for further references and some generalizations.

\section{The decomposability spectrum and $\pcf$ theory}\label{pcfth} 

\begin{theorem}\label{ldec} 
Suppose that $E$ is an ultrafilter over $J$, 
$(\lambda _j)_{j\in J}$ are regular cardinals and
$\cf \prod_E \lambda _j=\mu $.
If $D$ is a $|J|^+$-complete ultrafilter 
and $D$ is $\lambda_j$-decomposable for all
$j \in J$
then $D$ is $\mu $-decomposable.
\end{theorem}

\begin{proof}
Since 
$\cf \prod_E \lambda _j=\mu $,
there are functions $g_\alpha \in \prod \lambda_j$ ($\alpha\in \mu $)
such that $[g_\alpha ]_E$ ($\alpha\in \mu $) is a cofinal sequence
in $\prod_E \lambda_i$. Without loss of generality, we can assume
$[g_\alpha ]_E <_E [g_{\alpha'} ]_E $ whenever $\alpha < \alpha' < \mu $.

For every $j\in J$, let $f_j$ be a $\lambda_j$-decomposition of $D$, and 
suppose that
$D$ is over $I$.
Thus, for every $j\in J$ and for every $\beta $ in $\lambda _j$ we have
$\{i \in I| f_j(i)>\beta  \} \in D $. In particular, for every 
 $\alpha \in \mu $ and $j\in J$ we have $\{i \in I| f_j(i)>g_\alpha (j) \} \in D $.
Since $D$ is  $|J|^+$-complete, then
for every $\alpha \in \mu $, it happens that
 $X_\alpha =\bigcap_{j \in J} \{i \in I| f_j(i)>g_\alpha (j) \} \in D $.

Thus, if $i\in X_\alpha $ then $ f_j(i)>g_\alpha (j)$ for all $j \in J$.
Hence, if $i\in X_\alpha $, then $ [f_j(i)]_E>_E[g_\alpha]_E$.
If we put $Y_\alpha=\{i \in I| [f_j(i)]_E>_E[g_\alpha]_E\} $
then $Y_\alpha\supseteq X_\alpha \in D$, hence 
$Y_\alpha \in D$, for every $\alpha \in \mu $.

Since $[g_\alpha ]_E <_E [g_{\alpha'} ]_E $ whenever $\alpha < \alpha' < \mu$,
we have $Y_\alpha \supseteq Y_{\alpha'}$ whenever $\alpha < \alpha' < \mu $.
Moreover, $\bigcap_{\alpha \in \mu } Y_\alpha = \emptyset$, since
if, on the contrary, $i \in \bigcap_{\alpha \in \mu } Y_\alpha $, then
$[f_j(i)]_E>_E[g_\alpha]_E $ for all $\alpha \in \mu $, and this contradicts
the assumption that $[g_\alpha ]_E$ ($\alpha\in \mu $) is a cofinal sequence
in $\prod_E \lambda_j$.

Thus, we have found a sequence $ Y_\alpha $ ($\alpha\in \mu $)
of sets in $D$ such that 
$Y_\alpha \supseteq Y_{\alpha'}$ whenever $\alpha < \alpha' < \mu $,
and $\bigcap_{\alpha \in \mu } Y_\alpha = \emptyset$. This means
that $D$ is $\mu $-descendingly incomplete. Since $\mu$ is a 
regular cardinal by assumption  (being the cofinality of $\prod_E \lambda _i $)
we get that $D$ is $\mu $-decomposable.
\end{proof}

In the above theorem we are not necessarily assuming that
all the $ \lambda _i$'s are distinct. The particular case in which they are
all distinct can be, of course, restated in terms of $\pcf$ theory \cite{She}. 

\begin{corollary}\label{pcfl}
Suppose that $\aaa$ is a set of regular cardinals, $\mu \in \pcfa$, 
$D$ is an $|\aaa |^+$-complete ultrafilter, and $D$ is $\kappa $-decomposable for every $\kappa \in \aaa$.
Then $D$ is $\mu$-decomposable.
\end{corollary}

%\begin{proof}\end{proof}
%Recall \cite{She} that $\mu \in \pcfa$ means that. 

Notice that the condition $ |\aaa|  < \min \aaa $ follows from the 
hypotheses of Corollary \ref{pcfl}, since if the ultrafilter   
$D$ is $|\aaa |^+$-complete and $\kappa $-decomposable then necessarily
$ \kappa > |\aaa|$. 

%By applying $\pcf$-theory, 
%%  ??Theorem \ref{ldec} and 
%Corollary \ref{pcfl} 
%can be used to obtain partial solutions of some 
%problems from \cite{L1}.

We now show that, for $ \lambda $ singular, 
if $K_D \cap \lambda $ is an interval of regular cardinals 
cofinal in $ \lambda $, and $D$ is
 $ (\cf  \lambda )^+$-complete, then 
$\lambda^+ \in K_D$. In turn, 
$ \lambda ^+$-decomposability,
together with $ (\cf  \lambda )^+$-completeness, implies 
$\pp( \lambda )= \lambda ^+$. 
This means that, in some sense, the $\pcf$ theory at
$ \lambda $ has the simplest possible structure, and, as argued in 
\cite{She}, it is a version of the Generalized Continuum Hypothesis.
Moreover, $\pp( \lambda )= \lambda ^+$ and Corollary \ref{pcfl}
imply that we can equivalently suppose that there are arbitrarily large
regular $ \kappa < \lambda $ such that $D$ is $ \kappa $-decomposable,
as exemplified in Proposition \ref{solovaygenpp} below. 
A $\pcf$-free version of Theorem \ref{limit} will be given in 
Proposition \ref{solovaygen}.   

\begin{thm}\label{limit}
Suppose that $ \lambda  $ is a singular cardinal, and 
the ultrafilter $D$ is $ (\cf  \lambda )^+$-complete and 
$ \kappa $-decomposable for all sufficiently large regular $ \kappa < \lambda $.

Then $D$ is $\lambda^+$-decomposable, and moreover $\pp( \lambda )= \lambda ^+$.   
\end{thm}

\begin{proof}
By \cite[II, Theorem 1.5]{She} there is a strictly increasing sequence 
$ \lambda _ \alpha $ ($ \alpha  \in \cf \lambda $) of regular cardinals
 with $ \lambda _ \alpha  < \lambda $, and
 $\sup _{ \alpha  \in \cf \lambda }  \lambda _ \alpha = \lambda $, 
and there is a uniform ultrafilter $E$ on $\cf \lambda $ such that 
$ \cf \prod_E \lambda _ \alpha = \lambda ^+$. Actually, \cite[II, Theorem 1.5]{She}
obtains the above result  for the ideal $J ^{\bd} _{\cf \lambda } $ in place of the ultrafilter $E$; however, it is enough
to take as $E$ any ultrafilter extending the dual of  
 $J ^{\bd} _{\cf \lambda } $ 
(cf. e. g. the proof of \cite[Lemma 1.4]{BM}).
%%%Dire tutto questo nell'introduzione, visto che lo applico piu' volte?
Notice that the fact that  $J ^{\bd} _{\cf \lambda } $ is the ideal of sets bounded in $ \cf \lambda $ implies that $E$ is uniform over
$ \cf \lambda $.     

By assumption, 
there is $ \lambda' < \lambda $ such that 
 $D$ 
is $ \kappa $-decomposable
for every regular $ \kappa  $ with $ \lambda  ' \leq \kappa < \lambda $.
Since $E$ is uniform, the set
 $X= \{ \alpha \in \cf \lambda | \lambda ' \leq \lambda _ \alpha \} $
belongs to $  E$. Hence, if $E'$ is the restriction of $E$ to $X$,
the hypotheses of Theorem \ref{ldec} apply with $E'$
in place of $E$, and we get that $D$ is $ \lambda ^+$-decomposable.

We now show that if there is a $ \lambda ^+$-decomposable
and $ (\cf  \lambda )^+$-complete ultrafilter $D$ then
$\pp( \lambda )= \lambda ^+$.
%By a remark in the introduction,
Without loss of generality, we can suppose that
$D$ is uniform on $ \lambda ^+$.
Since $D$ is $ (\cf  \lambda )^+$-complete,
$D$ is not  $ \cf  \lambda $-decomposable, hence,
by Theorem \ref{cckp}, there is $ \lambda ' < \lambda $
such that $D$ is $(\lambda',\lambda^+ )$-regular.
%By a remark in the introduction,
Hence, $D$ is $ \kappa $-decomposable for all 
regular $ \kappa $'s with $ \lambda ' \leq \kappa < \lambda $.   

Suppose by contradiction that $\pp( \lambda )> \lambda ^+$.
By the ``No Hole Conclusion'' \cite[II, 2.3(1)]{She},   
$ \lambda^{++} \in \pcfa$, where $\aaa$ is
a set of regular cardinals cofinal in $ \lambda $, and 
$|\aaa|= \cf \lambda $ (again \cite{She} deals with an ideal,
but it is sufficient to extend the dual of this ideal to an ultrafilter,
which turns out to be uniform over $\cf \lambda $).
By considering, as above, a final segment of $ \aaa$, 
we get from Corollary \ref{pcfl} that $D$ is
$ \lambda ^{++}$-decomposable, but this is impossible, 
since $D$ is uniform over $ \lambda ^+$.    
\end{proof}

\begin{corollary} \label{cfomega}
If $ \lambda $ is a singular cardinal of cofinality $ \omega$, and
 the ultrafilter $D$ is
$ \kappa $-decomposable for all sufficiently large regular $ \kappa < \lambda $
then $D$ is either $ \lambda $-decomposable, or $ \lambda ^+$-decomposable.
\end{corollary}

\begin{proof}
If $D$ is $ \omega _1$-complete, then $D$ is $ \lambda ^+$-decomposable, by Theorem  \ref{limit}.

On the other side, if $D$ is not $ \omega _1$-complete, then
 $D$ is $ \omega $-decomposable, hence $ \lambda $-decomposable by
\cite[Proposition 1, and footnote on p. 461]{P1}. See also \cite{L1} 
for generalizations of results from \cite{P1}.   
\end{proof}

\begin{proposition} \label{solovaygenpp}
If $ \lambda  $ is a singular cardinal
and the ultrafilter $D$ is $ (\cf  \lambda )^+$-complete, then the following conditions are equivalent:

(a) $D$ is $ \kappa $-decomposable for all sufficiently large regular $ \kappa < \lambda $;

(b) $D$ is $\lambda^+$-decomposable;

(c) There is some $ \lambda' < \lambda $ such that 
$D$ is $(\lambda',\lambda^+ )$-regular;

(c$'$) There is some $ \lambda' < \lambda $ such that 
$D$ is $(\lambda',\lambda)$-regular;
 
(d)  $\pp( \lambda )= \lambda ^+$ and there
are arbitrarily large regular cardinals $ \kappa < \lambda $
such that $D'$ is $ \kappa $-decomposable.
%(e) $\pp( \lambda )= \lambda ^+n$ etc.???
\end{proposition}  

\begin{proof} 
(a) $ \Rightarrow $ (b) $ \Rightarrow $ (c)  and 
(b) $ \Rightarrow $ (d) are given
by the proof of Theorem \ref{limit}.

(d) $ \Rightarrow $ (b) follows from Corollary \ref{pcfl}, and 
 (c) $ \Rightarrow $ (c$'$) $ \Rightarrow $ (a) are trivial.
%follows from a remark in the introduction.
\end{proof}

\section{A generalization of Solovay's GCH result}\label{Sol} 
% fare un'unica sezione "further remarks?"

Some arguments from the proof of  Theorem  
\ref{limit} (with no use of $\pcf$ theory)
can be used to furnish an alternative proof of
 R. Solovay's 
GCH result.

\begin{proposition} \label{solovaygen}
Suppose that $ \lambda  $ is a singular cardinal, 
$ \cf \lambda < \kappa$, and 
$\nu ^{<\kappa} < \lambda $ for every $\nu< \lambda $.
If there exists an ultrafilter $D$ which is 
$\lambda^+$-decomposable and
$\kappa$-complete, 
 then $ (\lambda ^+)^{<\kappa}=\lambda ^+$.
\end{proposition}  

\begin{proof}
Without loss of generality,
%as in
we can suppose that $D$ is uniform over $\lambda ^+$.
Since $D$ is 
$\lambda^+$-decomposable
 and  not $ \cf  \lambda $-decomposable
(being $\kappa$-complete and $\kappa>\cf \lambda $), then,
by Theorem \ref{cckp}, 
$D$ is $(\nu, \lambda ^+)$-regular 
for some $ \nu < \lambda $.

Since $D$ is $\kappa$-complete then,
by a remark which is probably due to Solovay
(see 
%\cite[p. 74]{Ket1}, 
\cite{L1}%citare il posto
),
$D$ is $((\nu^{<\kappa})^+, (\lambda ^+)^{<\kappa})$-regular.
Since $\nu^{<\kappa}<\lambda $, and $D$ is uniform over
$\lambda^+$, this can happen only if $(\lambda ^+)^{<\kappa}=\lambda^+$.
% by a remark in the introd???
\end{proof}

\begin{theorem}\label{Solov}
\cite{So}
Suppose that $\kappa $ is $\mu $-strongly compact.

(a) If $\kappa \leq \nu \leq \mu $, and $ \nu $  is regular, then $\nu ^{< \kappa}=\nu  $.

(b) If $\kappa \leq \lambda  \leq \mu $, and $ \lambda $
is singular and strong limit, then $2^\lambda =\lambda ^+$.
\end{theorem}

\begin{proof}
That $\kappa $ is $\mu $-strongly compact implies
that there is a $ \kappa $-complete $( \kappa , \mu)$-regular ultrafilter.
Hence, for every regular $ \nu $ with $\kappa \leq \nu \leq \mu $
there is a $ \kappa $-complete $ \nu $-decomposable ultrafilter. 

(a) is now proved by induction on  $ \nu $.
Since $ \kappa $ is measurable, hence strongly inaccessible,
and because of standard cardinal arithmetic, the only non trivial case is
when $\nu$ is the successor of a singular cardinal $ \lambda $ of cofinality
$ < \kappa $. In this case, we can apply Proposition \ref{solovaygen}
because of the above remark.   

(b) In the case when $ \cf \lambda < \kappa $, case (b) follows
from (a) by standard cardinal arithmetic. The case   
$ \cf \lambda \geq \kappa $ can be obtained  as a consequence of 
Silver's Theorem (see \cite[p. 191]{KM}) from the case $ \cf \lambda < \kappa $, or, 
alternatively, using the appropriate arguments from \cite{So}.
\end{proof}

\section{Some problems}\label{prob} 

\begin{prob}\label{conject}
If we assume that $E$ is uniform, and that  $|J|< \lambda _i$ for all $i$,
can the 
hypothesis ``$D$ is 
$|J|^+$-complete'' in Theorem \ref{ldec}
 be weakened to ``$D$ is not $|J|$-decomposable''?
\end{prob}
 
The assumption ``$|J|< \lambda _i$ for all $i$'' in Problem
\ref{conject}  is necessary, otherwise there are 
easy counterexamples.
However,
the assumption that $D$ is 
$ (\cf \lambda )^+$-complete
 can be weakened to 
``$D$ is not $\cf \lambda $-decomposable'',
in Theorem  \ref{limit} and 
Proposition  \ref{solovaygenpp}.
%altri??

\begin{proposition} \label{solovaygenppy}
If $ \lambda  $ is a singular cardinal
and the ultrafilter $D$ is not $ \cf  \lambda $-decomposable,
then the following conditions are equivalent:

(a) $D$ is $ \kappa $-decomposable for all sufficiently large regular $ \kappa < \lambda $;

(b) $D$ is $\lambda^+$-decomposable;

(c) There is some $ \lambda' < \lambda $ such that 
$D$ is $(\lambda',\lambda^+ )$-regular;

(c$'$) There is some $ \lambda' < \lambda $ such that 
$D$ is $(\lambda',\lambda)$-regular;
 
(d)  $\pp( \lambda )= \lambda ^+$ and there
are arbitrarily large regular cardinals $ \kappa < \lambda $
such that $D'$ is $ \kappa $-decomposable.
%(e) $\pp( \lambda )= \lambda ^+n$ etc.???
\end{proposition}  

A major problem is what happens if 
in the hypothesis of Theorem  \ref{limit}, Corollary
 \ref{cfomega} and Proposition \ref{solovaygenpp}(a) we only suppose 
that there are arbitrarily large 
regular cardinals $ \kappa < \lambda $ such that  
 $D$ 
 is $ \kappa $-decomposable.
We have partial results showing that the statements still hold,
except possibly for very special situations.

We can ask even a subtler problem.
Does the first conclusion in Theorem  \ref{limit} hold when 
decomposability is replaced by regularity?
Namely, is the following true?

\begin{conj}\label{conjlimit}
Suppose that $ \lambda  $ is a singular cardinal, $\mu < \lambda $ and 
the ultrafilter $D$ is $ (\cf  \lambda )^+$-complete and 
$(\mu, \kappa )$-regular for all $ \kappa < \lambda$.

Then $D$ is $(\mu, \lambda ^+)$-regular.
\end{conj}

Again, we expect to get this, except perhaps for really
special situations, and it is likely that the assumption
that $D$ is $ (\cf  \lambda )^+$-complete can be weakened to 
``$D$ is not $\cf \lambda $-decomposable''.  

A positive solution to the above problems would furnish a solution to
a lot of problems raised in \cite{L1}.  
%Non ripeterlo diecimila volte!!!!


\begin{thebibliography}{AB} 

\bibitem[BM]{BM} M. Burke and M. Magidor, \emph{Shelah's ${\rm pcf}$ theory and its applications}. 
   Ann. Pure Appl. Logic {\bf 50} (1990), no. 3, 207--254. 

\bibitem[CC]{CC} G. V. Cudnovskii and D. V. Cudnovskii, 
\emph{Regular and descending 
incomplete ultrafilters} (English translation), Soviet Math. Dokl. 
{\bf 12} 901--905 (1971).

\bibitem[D]{D} H. D. Donder, \emph{Regularity of ultrafilters and the core model}, Israel 
J. Math. {\bf 63} 289--322 (1988).

\bibitem[KM]{KM} A. Kanamori and M. Magidor, \emph{The evolution of large cardinal axioms in 
Set Theory}, in \emph{Higher Set Theory} (G. H. M\"uller and D. S. Scott editors) 
99--275, Berlin (1978).

\bibitem[KP]{KP}  K. Kunen and K. L. Prikry, \emph{On descendingly incomplete ultrafilters}, 
J. Symbolic Logic {\bf 36} 650--652 (1971).

\bibitem[L]{L1} P. Lipparini, \emph{More on regular ultrafilters in ZFC}, to be revised, preliminary version available at the author's web page.

\bibitem[P]{P} K. Prikry, \emph{Changing measurable into accessible cardinal}, 
Dissertationes Mathematicae (Rozprawy Matematyczne) LXVIII (1970). 

\bibitem[P1]{P1} K. Prikry, \emph{On descendingly complete ultrafilters}, 
in \emph{Cambridge 
Summer School in Mathematical Logic} (A. R. D. Mathias and H. Rogers editors) 
459--488, Berlin (1973).   

\bibitem[Shr]{Shr}   M. Sheard,
\emph{Indecomposable ultrafilters over small large cardinals},
   J. Symb. Logic {\bf 48} 1000--1007 (1983).

\bibitem[She]{She} S. Shelah, \emph{Cardinal Arithmetics}, Oxford (1994).

\bibitem[Si]{Si}
 J. H. Silver,
 \emph{Indecomposable ultrafilters
   and $0^\#$}, in \emph{\it  Proceedings of the Tarski Symposium},  Proc. Sympos. 
Pure
   Math. XXV, Univ. Calif., Berkeley, Calif., 357--363 (1971).
 
\bibitem[So]{So} Solovay, \emph{Strongly compact cardinals and the GCH.} in 
   \emph{Proceedings of the Tarski Symposium} (Proc. Sympos. Pure Math., Vol.
   XXV, Univ. California, Berkeley, Calif., 1971), pp. 365--372.
   Amer. Math. Soc., Providence, R.I., 1974.


\end{thebibliography}
\end{document}